\documentstyle[12pt]{article}

\title{Anchored expansion and random walk}
\author{B\'alint Vir\'ag}
\date{January 19, 2000}
\newtheorem{thm}{Theorem}[section]
\newtheorem{pro}[thm]{Proposition}
\newtheorem{lem}[thm]{Lemma}
\newtheorem{cor}[thm]{Corollary}
\newenvironment{thm*}[1]{\par \trivlist
 \itemindent 0pt \item[\hskip\labelsep\bf Theorem #1]
 \it\ignorespaces}{\endtrivlist}
\newenvironment{eg}{\par \trivlist
 \itemindent 0pt \item[\hskip\labelsep\bf Example
 \refstepcounter{thm} \thethm]\ignorespaces}{\endtrivlist}
\newcommand{\qed}{\hfill\mbox{$\framebox(5,5)[]{}$}}
\newenvironment{proof}{\par \trivlist
 \itemindent\parindent \item[\hskip\labelsep\sc Proof.]
 \ignorespaces}{\qed\endtrivlist}

\newcommand{\Gr}{{\cal G}}
\newcommand{\walk}{$\{X_n\}$}

\newcommand{\eps}{\varepsilon}
\newcommand{\Var}{\mbox{\rm Var}}

\newcommand{\hf}{\frac{1}{2}}
\newcommand{\ev}{\mbox{\bf E}}
\newcommand{\pr}{\mbox{\bf P}}
\newcommand{\one}{\mbox{\bf 1}}
\newcommand{\as}{\mbox{\hspace{.3cm} a.s.}}
\newcommand{\dist}{\mbox{\rm dist}}
\newcommand{\is}{\Delta}
\newcommand{\wmax}{{w_0}}
\newcommand{\aec}{\mbox{\bf i}}
\newcommand{\dd}{\mbox{\boldmath $\delta$}}
\newcommand{\ctry}{{\Xi}}
\newcommand{\disti}{\dist_i}
\newcommand{\sfrac}[2]{\mbox{$\frac{#1}{#2}$}}
\newcommand{\lsp}{{\underline S}}
\newcommand{\rad}{r}

\begin{document}
\maketitle
\begin{abstract}
This paper studies anchored expansion, a non-uniform version of the
strong isoperimetric inequality. We show that every graph with
$i$-anchored expansion contains a subgraph with isoperimetric
(Cheeger) constant at least $i$. We prove a conjecture by Benjamini,
Lyons and Schramm (1999) that in such graphs the random walk escapes
with a positive lim inf speed. We also show that anchored expansion
implies a heat-kernel decay bound of order $\exp(-cn^{1/3})$.
\end{abstract}
\footnotetext[1]
{{Research  partially supported by
the NSF grant \#DMS-9803597 and the the Lo\`eve Fellowship}}
\footnotetext[2]
{{\it AMS} 1991 {\it subject classifications. } Primary 60J15, 60F15;
secondary 05C75.}

\section{Introduction}

Anchored expansion was introduced by Benjamini, Lyons, and Schramm
(1999) as a non-uniform version of the strong isoperimetric
inequality, after Thomassen (1992) used more general ``anchored''
isoperimetric inequalities to give sufficient conditions for the
transience of random walks on graphs. Let $G=(V,w)$ be an infinite
weighted graph, that is a countable set $V$ together with a
symmetric, nonnegative function $w((u,v))$, and let $E$ denote the
set of edges, that is unordered pairs $(u,v)$ where $w$ is positive.
Define the {\bf weight} $w(v)$ of a vertex $v$ as the sum of the
weights over the incident edges; we will assume that this is finite
for every vertex. Define the {\bf volume} $|\cdot|$ of an edge or
vertex set as the sum of the weights over the set, and define the
{\bf edge boundary} $\partial S$ of a vertex set $S$ as the set of
edges with one vertex inside $S$ and one outside.  The {\bf strong
isoperimetric inequality} with constant $i$, perhaps the simplest
isoperimetric inequality, states that
 \begin{equation} \label{ip}
 |\partial S| \ge i|S|
 \end{equation}
 for all vertex sets $S$, and the {\bf Cheeger constant} of a graph
is the supremum of the $i$ for which this inequality holds.  The
implications of a positive Cheeger constant are very strong; the
Markov kernel on a graph with positive Cheeger constant has spectral
radius less than 1 (Cheeger (1970), Dodziuk (1984), Mohar (1988);
these two conditions are in fact equivalent) and as a result, if the
graph does not grow faster than exponentially, the random walk
escapes at a positive speed, (i.e. linear rate). As is clear from
the definition, having a positive Cheeger constant is a rather
fragile property. Random perturbations, such as Bernoulli
percolation on an (unweighted) graph, even with a very high 
retention probability, or
a geometric stretching of edges, destroy it; see BLS (1999) for
results on stability of graph properties under random perturbations.

 A more stable condition, which we call {\bf $i$-anchored
expansion}, requires (\ref{ip}) to hold with possibly some
exceptions provided that every vertex is contained in only finitely
many connected exceptional $S$ (of course, there need not be a
uniform bound on the number of such exceptions). We call the
supremum of $i$ for which a graph $G$ has $i$-anchored expansion the
{\bf anchored expansion constant} $\aec (G)$; if this constant is
grater than 0, then we say that $G$ has the {\bf anchored expansion
property}. If $G$ is connected with edge weights bounded from below
(in particular, if $G$ is a connected, unweighted graph), then this
is equivalent to a version of the original definition,
 $$
 \aec(G) = \liminf_{S \ni v} \frac{|\partial S|}{|S|},
 $$
 here $S$ ranges over all connected vertex sets containing a fixed
vertex $v$;  the lim inf is then applied to the set of values
obtained. This definition does not depend on the choice of the fixed
vertex, and it explains the name ``anchored expansion''. The only
difference from the original definition (given for the unweighted
case) is that here the volume of a vertex set is defined as the sum
of degrees rather than the number of vertices.

 Benjamini, Lyons and Schramm (1999) conjectured that in an
unweighted graph with bounded degree and anchored expansion, the
random walk has positive lim inf speed with positive probability.
The main goal of this paper is to prove this conjecture in a
slightly stronger form. Say a weighted graph has {\bf
$\wmax$-bounded geometry} if all positive edge weights are at least
1 and all vertex weights are at most $\wmax$, and let $|X_n|$
denote the graph distance between the random walker and the initial
vertex at time $n$.
 \def\mainth{There exists $c>0$ so that the random walk on a
weighted graph with $\wmax$-bounded geometry satisfies $\liminf
|X_n|/n \ge c\aec(G)^7\wmax^{-3}$ a.s.}
 \begin{thm} \label{main th} \mainth \end{thm}
 This theorem gives a geometric explanation for positive speed in
certain graphs, such as infinite components of $p$-Bernoulli
percolation on graphs with positive Cheeger constant for high $p$,
geometric edge-stretchings of such graphs, or supercritical
Galton-Watson trees. For these graphs, Chen and Peres (1998),
inspired by questions of Benjamini, Lyons, and Schramm (1999),
proved the anchored expansion property. It is an open question (see
H\"aggstr\"om, Schonmann and Steif (1999)) whether infinite clusters
of Bernoulli percolation on a transitive graph have the anchored
expansion property. The same authors prove that if $G$ is a
transitive graph and there exists an automorphism-invariant
percolation on $G$ where all infinite components have the anchored
expansion property, then $G$ has positive Cheeger constant.

 Since an exponential heat kernel bound is equivalent to having a
positive Cheeger constant, one cannot hope that anchored expansion
would imply such a strong bound. The following theorem gives a
sub-exponential bound, which is strongest in the sense that the
$n^{1/3}$ in the exponent cannot be improved.

 \def\heatth{Let $G$ be a weighted graph with the anchored expansion
property and $\wmax$-bounded geometry. Let
 $\alpha:=\aec (G)^2(\wmax^2/2)^{-{1/3}}/9$.
 For every vertex $x$ there is an $N$ so that
 $$
 p^n(x,y) < e^{-\alpha n^{1/3}}
 \ \ \ \mbox{for all\ } n>N,\ y\in V.
 $$}
 \begin{thm} \label{heat th} \heatth \end{thm}
 Varopoulos (1991) showed that such a bound holds in Cayley graphs
of any finitely generated group of exponential growth (see Hebisch
and Saloff-Coste (1993) for a more general statement and a simpler
proof). Such a group can be amenable (e.g. the lamplighter group
$G_1$), in which case it is an example where such decay holds but
anchored expansion does not.

Section \ref{geometry s} examines the geometry of graphs with
$i$-anchored expansion, and shows that they are built from a graph
with Cheeger constant at least $i$ and ``islands'', each of which is
finite but whose size is not necessarily bounded by a constant (a
binary tree with ``pipes'' of increasing length attached to a scarce
set of vertices is a graph with anchored expansion and unbounded
``islands'').  Section \ref{walk s} proves some properties of random
walk on such graphs. Using these results, in Section \ref{speed s}
we prove that the random walk has positive speed, and in Section
\ref{heat s} we establish the heat kernel bound.

\section{Notation}

The concept of {\bf the random walk on a weighted graph} $G=(V,w)$
is just a geometric representation of a countable, reversible Markov
chain with transition probabilities
 $$
 p(u,v)=w((u,v))/w(u).
 $$ 
 Conversely, if we have a reversible Markov chain with transition
probabilities $p(u,v)$ and stationary measure $w(v)$ we get a
weighted graph by the above formula. We will use the notation $P$
for the transition probability matrix of the walk. We will usually
denote the walk itself by $\{X_n\}$, and $\pr_v$, $\ev_v$ will mean
probability and expectation with respect to the walk started at
vertex $v$.

We will consider the Hilbert space $L^2(V,w)$ of functions, equipped
with the inner product and norm
 $$
 (f,g)=\sum_{v\in V}f(v)g(v)w(v), \ \ \ \
 \|f\|=(f,f)^\hf.
 $$ 
 For an operator $P$ on $L^2(V,w)$, we will use the norm $\|P\|=
\sup \|Pf\|/\|f\|$. The {\bf Markov kernel} $P$ of a random walk on
a weigthed graph is the operator on $L^2(V,w)$ defined by $(Pf)(v):=
\ev_v f(X_1)$.

We will use the standard notation $\dd_x$, for unit mass at $x$
(formally, $\dd_x$ is an element of the dual of $L^2(V,w)$), and
$\one_S$ for the indicator (characteristic) function of the set $S$.  
The volume of a vertex set $S$ can be written as $|S| = \sum_{v\in
S} w(v) = \|\one_S\|^2$. The {\bf number of vertices} of $S$ will be
denoted $\#S$.

The {\bf induced subgraph} of $S\subset V$ is the graph $G$ with
vertex set $S$ and edge set given by $\{(u,v)\in E: u,v \in S\}$,
and we call $S$ {\bf connected} if its induced subgraph is
connected. We will often write $G\setminus S$ for the induced
subgraph of $V\setminus S$. The {\bf inner vertex boundary} of $S$
is the set of vertices in $S$ with a neighbor outside $S$, and the
{\bf outer vertex boundary} of $S$ the set of vertices outside $S$
with a neighbor in $S$.

By a {\bf path} of {\bf length} $n$ we mean a subgraph with vertex
set $(v_0,\ldots, v_n)$ and edge set consisting of edges between
consecutive $v$-s.  The {\bf graph distance} between two vertices in
$G$ is the length of the shortest path with endpoints given by the
two vertices. The notation $|v|$ for a vertex denotes the graph
distance between $v$ and some declared fixed vertex; this vertex
will usually be the starting point of the random walk we study.  
The {\bf lim inf speed} of the walk is defined as $\liminf |X_n|/n$.

We will use $c,\ c_2,\ c_3$ to denote constants whose values might
change from one expression to another.

\section{Geometry of graphs with anchored expansion }
\label{geometry s}

 Let $G=(V,w)$ be a weighted graph. For $i<1$, define the {\bf
$i$-isolation} of a finite-volume vertex set $S$ of $G$ by
 $$
 \is_i S = i|S|-|\partial S|.
 $$
 A vertex set $S$ with positive $i$-isolation will be called {\bf
$i$-isolated}. We will call a vertex set $S$ satisfying $\is_i
S>\is_i A$ for every subset $A\not=S$ of $S$ an {\bf $i$-isolated
core}. Since $A$ can be the empty set, an isolated core must be
either empty or isolated. A nice property of $i$-isolated cores is
given in the following lemma.

 \begin{lem}\label{isgrow}
 Let $A$ be a vertex set and let $S$ be an $i$-isolated core which is
not a subset of $A$. Then $\is_i (A \cup S) > \is_i A$.
 \end{lem}
 \begin{proof}
 Note that if $B$ and $C$ are disjoint vertex sets then
 \begin{equation}\label{union eq}
 \is_i (B\cup C)=\is_i B+\is_i C + 2|\partial B \cap \partial C|.
 \end{equation}
 The factor $2$ in the above expression appears since common
boundary edges of $B$ and $C$ are not boundary edges of their union.
Then
 $
 \is_i (A \cup S) =  \is_i (A \setminus S)+ \is_i S + 2|
\partial (A\setminus S)\cap \partial S|
 $.
 The hypothesis can be used to bound the second term. The last term
equals twice the total weight on edges with one endpoint in
$A\setminus S$ and the other endpoint in $S$; this does not increase
if we require the latter endpoint to be in a subset of $S$.
Therefore
 $$
  \is_i (A\cup S) > \is_i (A \setminus S)+ \is_i(A \cap S)  + 2|
\partial (A\setminus S)\cap \partial (A\cap S)| = \is_i A. \ \ \qed
 $$ \def\qed{} \end{proof}
 \begin{cor}\label{Ai}
 The union of finitely many $i$-isolated cores is an $i$-isolated core.
 \end{cor}
 \begin{proof}
 It suffices to prove this for two $i$-isolated cores $S,S'$. Let
$A\subset S\cup S'$. Two applications of the lemma imply $\is_i A
\le \is_i (A \cup S) \le \is_i (A\cup S \cup S')$, and one
inequality is strict unless $A=S\cup S'$.
 \end{proof}

Let $ A_i$ denote the union of all $i$-isolated cores in $G$. It
follows from the definitions and Corollary \ref{Ai} that if $G$ has
$i$-anchored expansion then every connected component of $ A_i$ is
a {\it finite} union of isolated cores, hence an $i$-isolated core.

 If $G$ has $i$-anchored expansion, then the set $ A_i$ has the
remarkable property that $G\setminus  A_i$ is a graph with Cheeger
constant at least $i$. Indeed, let $S$ be a finite subset of
$V\setminus  A_i$, and let $C$ be a (possibly empty) $i$-isolated
core containing all vertices adjacent to $S$ in $ A_i$. A (possibly
empty) minimal subset $B$ of $C \cup S$ satisfying $\is_i B \ge
\is_i (C \cup S)$ must be an $i$-isolated core, hence
$B\subset A_i$, and thus $B\subseteq C$. Since $C$ is an
$i$-isolated core, we get
 \begin{equation} \label{island touch}
 \is_i (C \cup S) \le \is_i B \le \is_i C.
 \end{equation}
 Let $\is_i^{G\setminus  A_i}$ denote $i$-isolation of vertex sets
in the graph $G\setminus  A_i$. When $ A_i$ is removed from $G$, the
volumes of both $S$ and $\partial S$ decrease by $|\partial C \cap
\partial S|$. Thus we get
 $
 \is_i^{G\setminus  A_i}S=\is_i S + (1-i)|\partial C \cap \partial
S|
 $.
 Expressing $\is_i S$ by (\ref{union eq}) gives
 $$
 \is_i^{G\setminus  A_i}S
 = \is_i (C\cup S) - \is_i C - 2|\partial C \cap \partial S|
 + (1-i)|\partial C \cap \partial S|.
 $$
 This is at most $0$ by (\ref{island touch}), and we get the
required isoperimetric inequality. Thus we have shown

 \begin{pro}
 Every graph with $i$-anchored expansion contains a subgraph with
Cheeger constant at least $i$.
 \end{pro}

Note that $G \setminus  A_i$ is an isomorphism-invariant function of
the graph $G$. If, for example, $G'$ is a transitive graph and
$G\subset G'$ is a random subgraph whose law is invariant under a
group of automorphisms of $G'$, then the law of $G\setminus  A_i$ is
also invariant under this group.

We will call the connected components of $ A_i$ ($i$-){\bf islands},
and $G\setminus  A_i$ the {\bf oceans} (plural since $G\setminus
 A_i$ is not always connected). If $i'<i$, then we have $\is_{i} S =
(i-i')|S| + \is_{i'} S$, so $i'$-isolated sets are also
$i$-isolated, and if $A \subset S$ and $\is_{i'} A < \is_{i'} S$,
then $\is_i A < \is_i S$. In particular, $i'$-isolated cores are
$i$-isolated cores as well, giving $A_{i'} \subset  A_i$. Thus
decreasing $i$ has the effect of global warming: it raises the level
of the oceans. The following lemma gives an upper bound on how much
the level needs to be raised to sink certain islands.

 \begin{lem}\label{sink}
 Let $G$ be a connected graph with $i$-anchored expansion and edge
weights bounded below by 1. Let $S$ be a union of islands, each
having volume at most $i'^{-1}$ for some $i'>0$. Then $S\subset
V\setminus A_{i'}$.
 \end{lem}

\begin{proof}
 Any $i$-island in $S$ has boundary volume at least $1$ and positive
$i$-isolation, thus volume greater than $i^{-1}$. This gives $i'<i$.
Similarly, each $i'$-island has volume greater than $i'^{-1}$, so no
$i'$-island is a subset of an island in $S$. But since
$A_{i'}\subset A_{i}$, this implies that $A_{i'}\cap S=\emptyset$.
 \end{proof}

 Let $N_m$ denote the $m$th positive time the walk is in
$G\setminus A_i$. The process $\{X_{N_m}\}$ (often called the
induced Markov chain on $G\setminus  A_i$) is also a reversible
Markov chain, that is a random walk on a weighted graph $G_i$. If
$G$ is connected, then so is $G_i$. The vertex set of $G_i$ is given
by $V\setminus  A_i$, and its edge weight function satisfies
 $$
 w_i((u,v))=w(u)\pr_u(X_{N_1}=v).
 $$
 Clearly, for $u,v \in V\setminus A_i$ we have $w_i((u,v))=w((u,v))$
unless $u$ and $v$ are both on the outer vertex boundary of the same
island in $G$. It is also clear that for $v\in V(G_i)$, we have
$w_i(v)=w(v)$. The reversibility of the walk on $G$ implies that
$w_i$ is a symmetric function on the edges.

 The graph $G_i$ has the same vertex set as $G\setminus  A_i$, but
its edge and vertex weights are greater or equal.  We now show that
$G_i$ also has Cheeger constant at least $i$. To see this, let $S$
be a finite subset of $V\setminus A_i$ and follow the argument for
$G\setminus  A_i$ to get (\ref{island touch}). Let superscript $G_i$
on volume or $i$-isolation denote these quantities measured in the
base graph $G_i$.  Since we have $|S|^{G_i}=|S|$ and $|\partial
S|^{G_i} \ge |\partial S| - |\partial C \cap \partial S|$, it
follows by (\ref{union eq}) and (\ref{island touch}) that
 $$
 \is_i^{G_i} S \le \is_i S + |\partial C \cap \partial S|
 = \is_i (C \cup S) - \is_i C - |\partial C \cap \partial S| \le 0.
 $$

The upcoming analysis of random walks will need a rigorous
formulation of the idea that large islands cannot be very close to
each other. One could expect islands to have a neighborhood, whose
radius depends on the size of the island, within which there are no
other islands; or if this cannot be achieved, then at least one
could group nearby islands together to get such a configuration.
This is too optimistic as said, but Proposition \ref{countries} has
a similar decomposition, for which we first have to introduce some
tools.

A {\bf bridge structure} interconnecting a vertex set $S\subset
V(G)$ is a set of vertices $B$ so that $B\cup S$ is a connected 
vertex set. A
{\bf bridge} connecting two vertex sets $S_1,S_2\subset V(G)$ is a
vertex set $B$ so that the vertex set $B \cup S_1 \cup S_2$ has a
connected component intersecting both $S_1$ and $S_2$. Define the
{\bf $i$-length} of a bridge $B$ as the number of its vertices in
the ocean, $\#(B\setminus A_i)$. For a vertex set $S$ and a vertex
$v\not\in S$ let $\disti(v,S)$ equal $1$ plus the $i$-length of the
shortest bridge connecting $\{v\}$ and $S$; for $v\in S$, let
$\disti(v,S):=0$.

 \begin{lem}\label{connect} Let $G$ be a weighted graph with
$\wmax$-bounded geometry and $\aec(G)>i>0$ for some $i$. Let ${\cal
R}$ be a set whose elements are unions of $i$-islands, and let $v$
be a vertex. Suppose that for each $R \in {\cal R}$, there exists a
bridge structure $B$ which interconnects $R\cup\{v\}$ and satisfies
 \begin{equation} \label{connect eq}
 \wmax\#(B\cup\{v\}\setminus A_i)/|R| \le \aec(G)-i.
 \end{equation}
 Then ${\cal R}$ is finite.
 \end{lem}
 \begin{proof}
 For $R\in \cal R$, let $S$ denote $(B \cup \{v\})\setminus A_i$,
and let $A$ denote the union of islands intersecting $B\cup R\cup
\{v\}$. Then
 $$
 |\partial (A\cup S)|\le |\partial A|+|\partial S| \le i|A| + |S|.
 $$ 
 The bound on the first term of the second inequality holds since
$A$ is a union of islands. By (\ref{connect eq}) we have $|S|/|A|
\le \aec(G)-i.$ Therefore, using that $A$ and
$S$ are disjoint, we get
 $$
 \frac{|\partial (A\cup S)|}{|A\cup S|} 
 \le \frac{i|A|+|S|}
          {|A|+|S|}
 =   \frac{i+ |S|/|A|}
          {1+ |S|/|A|}
 \le \frac{\aec(G)}{1+\aec(G)-i} <\aec(G).
 $$
 By the anchored expansion property there are only finitely many
such sets $A \cup S$ containing $v$. The lemma follows. 
 \end{proof}

 \begin{pro}\label{countries}
 Let $G$ be a graph with $\aec(G)>0$ and $\wmax$-bounded geometry.
Let
 $$
 0<i\le \sfrac{2}{3}\aec(G), \ \
 \rad(\ell):= a2^\ell/\ell^2, \ \
 a:=\sfrac{3}{2}\pi^{-2}i\wmax^{-1}.
 $$
 For each positive integer $\ell $ there is a (possibly empty)
collection $\ctry_\ell$ of vertex sets $C$, which we call {\bf level
$\ell$ countries}, so that the following hold:
 \list{}{\itemsep 0in \topsep 1ex}
 \item[-]For each $\ell$ and $C\in \ctry_\ell$, the set $C\cap A_i$
is a union of $i$-islands, and is called the {\bf land} of the country $C$.
Its volume satisfies $|C\cap A_i|\in [2^{\ell-1},2^\ell)$.
 \item[-]For each $\ell$ and $C\in \ctry_\ell$, $C\setminus A_i=\{v \in V(G_i): \disti(v,C\cap A_i) \le
\rad(\ell)\}$, and this set is called the {\bf waters} of the country $C$.
 \item[-]Any two countries at the same level are disjoint.
 \item[-]Every $i$-island is a subset of some country.
 \item[-]Each vertex of $G$ is contained in at most finitely many
countries.
 \endlist
 \end{pro}
 \begin{proof}
 We start by constructing {\bf regions} $R$, which are islands or
unions of islands, together with bridge structures $B(R)$ connecting
these islands if they are disjoint.  First, for each $\ell\ge 1$,
label each $i$-island with volume in $[2^{\ell-1},2^\ell)$ as a
level $\ell$ region, and for these regions $R$, set
$B(R)=\emptyset$.

Define the {\bf waters} of a level $\ell$ region $R$ as $\{v \in
V(G_i): \disti(v, R) \le \rad(\ell)\}$. Then, for $\ell=1,2,\ldots$
(in this order), consider a maximal matching of pairs of level
$\ell$ regions whose waters intersect, and label the union $R$ of
each matched pair $(R_1,R_2)$ of regions a level $\ell+1$ region.
Set $B(R)$ to be the union of $B(R_1)$, $B(R_2)$ and a shortest
(minimal length) bridge connecting $R_1$ and $R_2$.

For every level $\ell$ region $R$ and vertex $v$ in $R$ or in its
waters, consider the bridge structure $B(v,R)$ given by the union of
$B(R)$ and a shortest bridge connecting $\{v\}$ and $R$. We have
 $$
 \#(B(v,R)\cup \{v\}\setminus A_i)
 \le \rad(\ell)
 + \sum_{n=1}^{\ell-1} 2^{\ell-1-n} \cdot 2 \rad(n)
 < a 2^{\ell} \sum_{n=1}^\infty n^{-2}
 = {i \over \wmax} 2^{\ell-2}.
 $$
 In the second expression the first term is an upper bound on the
length of a shortest bridge connecting $\{v\}$ and $R$ plus $1$.  
The first factor in the sum is an upper bound on the number of pairs
of level $n$ regions contained in $R$; the second factor is an upper
bound of the length of a shortest bridge connecting such a pair.

 Since then $\wmax\#(B(v,R)\cup\{v\}\setminus A_i)/|R|<i/2\le
\aec(G)-i$, it follows by Lemma \ref{connect} that each vertex is
contained in only finitely many regions or their waters. Therefore,
every island is contained in a maximal region, that is a region
which is not contained in any other regions. Call the union of a
maximal region and its waters a country of level of the
maximal region. Call the region itself the land of the country. This
construction clearly satisfies the properties claimed in the
proposition.
 \end{proof}

\section{Random walk and anchored expansion} \label{walk s}

 Let $\{X_n\}$ be the random walk on a graph $G$ with $i$-anchored
expansion and $\wmax$-bounded geometry. Our strategy for the
analysis of this walk will be to handle the time spent in the oceans
and in the islands separately. Let $N_m$ be the $m$-th positive time
when $\{X_n\}$ visits a vertex in $G\setminus  A_i$. We have seen
that $\{X_{N_m}\}$ is the random walk on the graph $G_i$, which has
Cheeger constant at least $i$. First, we show that
 \begin{equation}\label{speedi}
 \liminf |X_{N_m}|/m \ge {|\log (1-i^2)| \over \log \wmax}
 > {i^2 \over \log\wmax} \as
 \end{equation}
 For this, we first quote a version of the classical result of
Cheeger (1970), Dodziuk (1984) and Mohar (1988), to be found, for
example, in Lyons and Peres (1998).

 \begin{pro}\label{Cheeger}
 Let $P$ be the Markov kernel of the random walk on
a weighted graph with Cheeger constant at least $i$. Then $\|P\| \le
(1-i^2)^{1/2} \le (1-i^2/2)$.
 \end{pro}

This, together with the following lemma implies (\ref{speedi}).

 \begin{lem}
 Let $G$ be a weighted graph with $\|P\|<1$, and let $f$ be a
nonnegative vertex function so that
 $$
 g := \limsup |f^{-1}([0,n])|^{1/n} < \infty.
 $$
 Then $\liminf f(X_n)/n \ge -2\log \|P\|/ \log g$ a.s.
 \end{lem}

 If we set $f(v):=|v|$ in $G$ (this might be different from $|v|$
measured in $G_i$), then the bounded geometry property implies that
$g\le \wmax$, and the lemma applied to the walk ${X_{N_m}}$ on
$G_i$ implies (\ref{speedi}).  In a similar fashion, we get the
bound $i^2/\log g$ for the lim inf speed of random walks in graphs
with Cheeger constant at least $i$ and exponential growth rate at
most $g$.

 \begin{proof}
 For a small $\eps>0$, let $a := -2 \log \|P\|/\log(g+\eps)-\eps$.
We have
 \begin{eqnarray*}
 \pr_x[f(X_n) \le an]
 &=& \dd_x P^n \one_{f^{-1}([0,an])}
 = w(x)^{-1} (\one_x,P^n \one_{f^{-1}([0,an])}) \\
 &\le& w(x)^{-1} \|P^n\| \|\one_{f^{-1}([0,an])}\|
 \end{eqnarray*}
 For sufficiently large $n$, this is bounded above by
$w(x)^{-1}\|P\|^n(g+\eps)^{an/2}$, which is summable, so $f(X_n)>an$
eventually a.s. \end{proof}

 Our next goal is to bound the time spent in vertex sets, in
particular, islands.

 \begin{lem}\label{time}
 Let $G$ be a graph with $i$-anchored expansion, let $S$ be a vertex
set, and suppose $i'\le i$ is a constant so that $S$ is contained in
$G_{i'}$. Let $n$ be an integer, $x\in V\setminus A_i$ be a vertex
with $\disti(x,S)\ge n+1$, and let $T$ be the time the random walk
on $G$ spends in $S$. Then we have
 \begin{eqnarray*}
 \pr_x(T>0) &\le& 2w(x)^{-\hf} \ (1-i^2)^{n\over 2} \ i^ {-2} \
 |\partial S|^{\hf}, \\
 \ev_xT     &\le& 2w(x)^{-\hf} \ (1-i^2)^{n\over 2} \ i'^{-2} \
 |S|^{\hf}, \\
 \ev_xT^2   &\le& 8w(x)^{-\hf} \ (1-i^2)^{n\over 2} \ i'^{-4} \
 |S|^{\hf}.
 \end{eqnarray*}
 For an arbitrary vertex $x$, these bounds hold with $n=0$. If all
edge weights are at least $1$ and $S$ is a union of islands, then we
can use $i':=|S|^{-1}$.
 \end{lem}
 \begin{proof}
 The quantities $T$, $w(x)$, $|\partial S|$, $|S|$ do not change if
they are considered (for the walk) in the graph $G_{i'}$ instead of
the graph $G$, so we will do this.

Denote $P_i,\ \Gr_i,\ P_{i'},\ \Gr_{i'}$, the Markov and Green
kernels of the walks on $G_i$, and $G_{i'}$, respectively. Recall
that $\Gr_i= \sum_{m=0}^\infty P_i^m$, so we have
 $$
 \|\Gr_i\| \le \sum_{m=0}^{\infty} \|P_i\|^m = {1 \over 1- \|P_i\|},
 $$
 and so from Proposition \ref{Cheeger} we get
 \begin{equation}\label{norm bound}
 \|P_i\|\le (1-i^2)^\hf\le (1-i^2/2), \ \ \ \|\Gr_i\| \le 2i^{-2},
 \end{equation}
 and these inequalities also hold with $i$ replaced by $i'$
everywhere. For the walk on $G_i$ started at $y\in V\setminus A_i$,
the probability of moving into $S$ from the outside in one step is
given by the function $f(y)$ which equals $\dd_yP_{i'} \one_S$
outside $S$, and 0 in $S$.
 Thus the chance of moving into $S$ from the outside after $m$ steps
in $V\setminus A_i$ is given by $\dd_yP_i^mf$, and therefore
 $$
 \pr_x(T>0)
 \le \sum_{m=n}^{\infty} \dd_x P_i^m f.
 $$
 In Green kernel notation, this can be written as an inner product
 $$
 w(x)^{-1}(\one_x, P_i^n \Gr_i f)
 \le w(x)^{-1} \|\one_x\| \cdot \|P_i\|^n \cdot \|\Gr_i\| \cdot
\|f\|.
 $$
 The norms are all $L^2(V\setminus A_i, w)$, and the last inequality
follows from the Schwarz inequality and the norm bounds. Since $f\le
1$, the last norm is bounded above by $(f, 1)^\hf$, which equals
$|\partial S|^\hf$. The first claim of the lemma now follows from
(\ref{norm bound}).

For the expected value, write
 $$
 \ev_x T = \sum_{m=0}^{\infty} \dd_x P^n_iP_{i'}^m \one_S
 = w(x)^{-1}(\one_x, P^n_i\Gr_{i'} \one_S).
 $$
 Since $\|\one_S\|=|S|^\hf$, the norm bound on the last formula and
(\ref{norm bound}) give the second claim of the lemma.

 Finally, denote \walk\ the random walk on
$G_{i'}$. Then
 $$ \ev_x T^2
 = \ev_x \sum_{s,t > n \atop y,z\in S}\one(X_s=y)\one(X_t =z),
 $$
 and summing twice on or under the diagonal and extending the
range of $y$ gives the upper bound
 $$
 2 \ev_x \sum_{s>n \atop d
 \ge 0}\sum_{y\in V}\one(X_s=y)\one(X_{s+d} \in S).
 $$
 By the Markov property this equals
 $$
 2 \sum_{s>n \atop d \ge 0}\sum_{y\in V}\pr_x(X_s=y)\pr_y(X_{d}
\in S)
 = 2 \cdot \dd_x P_i^n \Gr_{i'} \Gr_{i'} \one_S.
 $$
 The third claim of the lemma follows if we write this as an inner
product and use norm bounds, as before. Omitting the estimates for
the first $n$ steps gives the proof for general $x$. Lemma
\ref{sink} implies that we can use $i':=|S|^{-1}$.
 \end{proof}

The anchored expansion property suggests that large islands cannot
be very frequent. The following lemma proves such a statement from
the point of view of the random walk. It uses the hypotheses and the
resulting decomposition of Proposition \ref{countries}.

 \begin{lem}
 Consider countries $C$ whose land is visited by time $N_m$, and let
$M_m$ be the volume of the largest such land. Then we have
 $$
 \lim\sup {M_m \over \log m (\log \log m)^2} < c< \infty \ \ \ \
\as
 $$
 \end{lem}
 \begin{proof}
 For a positive $b$, let $g(n):=b \log n (\log \log n)^2$, and let
${\cal A}_\ell$ be the event that the land of a country of level
$\ell$ is visited by time $g^{-1}(2^\ell)$. It suffices to prove
that only finitely many of these events happen, which will follow if
$\pr {\cal A}_\ell\le 2^{-\ell}$ for every large $\ell$. Consider
$\ell$ so large that the starting point of the walk is not contained
in any level $\ell$ country, and $\rad(\ell)\ge 1$. Let $T_C$ denote
the first hitting time of a country $C$. Let ${\cal A}_{\ell,C}$
denote the event that the land of the country $C$ is visited by time
$g^{-1}(2^\ell)$, and let ${\cal A}_C$ denote the event that the
land of the country $C$ is ever visited. The event ${\cal
A}_{\ell,C}$ implies ${\cal A}_C$ and $T_C\le g^{-1}(2^\ell)$, and
therefore
 $$
\pr {\cal A}_{\ell,C} \le \sum_{t=1}^{g^{-1}(2^\ell)}
\pr({\cal A}_C|T_C=t)\pr(T_C=t).
 $$ 
 Summing over level $\ell$ countries we get
 \begin{eqnarray}\nonumber
 \pr {\cal A}_\ell &\le&
 \sum_{C\in \ctry_\ell}\sum_{t=1}^{g^{-1}(2^\ell)} \pr({\cal A}_C|T_C=t)
 \pr (T_C=t)
 \\ &\le& \label{tudor} \sup_{C\in \ctry_\ell \atop 1\le t\le g^{-1}(2^\ell)}
 \pr({\cal A}_C|T_C=t)
 \sum_{t=1}^{g^{-1}(2^\ell)}\sum_{C\in \ctry_\ell} \pr (T_C=t).
\end{eqnarray}
 For fixed $t$, the events in the inner summand of (\ref{tudor}) are
disjoint, so the second factor is bounded above by $g^{-1}(2^\ell)$.
If $C$ is a level $\ell$ country with land $S$, then
 $$
 \disti({X_{T_C}}, S)=\lfloor \rad(\ell)\rfloor.
 $$
 Therefore by the Strong Markov Property and Lemma \ref{time},
$\pr({\cal A}_C|T_C=t)$ is not more than
 $$
 2w(x)^{-\hf}(1-i^2)^{(\rad(\ell)-2)/2}i^{-2}|\partial S|^\hf
 \le c' \exp(- c'_1 2^\ell / \ell^2) 2^{\ell/2}
 \le c 2^{-\ell} \exp(- c_1 2^\ell / \ell^2).
 $$

Then by (\ref{tudor}),
 $
 \pr {\cal A}_\ell
 \le g^{-1}(2^\ell)  c 2^{-\ell} \exp(- c_1 2^\ell / \ell^2)
 $
 and it suffices to prove that
 $$
 g^{-1}(2^\ell) \le c^{-1} \exp(c_1 2^\ell / \ell^2).
 $$
 We apply $g$ to both sides and use its monotonicity to transform
the above to
 $$
 2^\ell \le b (c + c_1 {2^\ell \over \ell^2}
 (\ell\log 2 - 2 \log \ell)^2 ).
 $$
 This certainly holds for all large $\ell$ if $b$ is large.
 \end{proof}

The following corollary will be used in a later section. It implies
that from the point of view of speed, distance can be measured while
walking on water.
 \begin{cor}
 \label{Hmcor}
 Set
 \begin{equation}\label{Hm}
 H_m:=\inf_{N_{m-1} < n \le N_m} |X_n|.
 \end{equation}
 Then we have
 $
 \lim (H_m/|X_{N_m}|) = 1
 $ a.s.
 \end{cor}
 \begin{proof} Between times $N_{m-1}$ and $N_m$ the walker is on an
island with diameter bounded above by the volume $M_m$ of the
largest land visited by time $N_m$. Thus we have $|X_{N_m}|-M_m-1
\le H_m \le |X_{N_m}|.$ Dividing by $|X_{N_m}|$, and using the lemma
together with (\ref{speedi})  proves the corollary.
 \end{proof}

\section{Lower bound on the speed}
 \label{speed s}

This section contains the proof of Theorem \ref{main th}. We also
give some counterexamples indicating why the bounded geometry
condition is important.
 \begin{thm*}{\ref{main th}} \mainth
 \end{thm*}
 \begin{proof}
 Let $G$ be a graph with the anchored expansion property and
$w_0$-bounded geometry, and consider the construction of countries
from Proposition \ref{countries}. Using the notation of the previous
section, we can decompose the inverse lim inf speed $\lsp^{-1}$ as
 \begin{eqnarray} \nonumber
 \limsup n/|X_n| &=&
 \limsup_m \sup_{N_{m-1} < n \le N_{m}} (n/|X_n|)
 \\ &\le& \limsup_m (N_m /H_m) = \limsup (N_m/|X_{N_m}|).
\label{sp1}
 \end{eqnarray}
 $H_m$ in the above expression is defined in (\ref{Hm}), and the
last equality follows from Corollary \ref{Hmcor}.
 Let $K_m:=N_m-m$ denote the time spent in the islands up to time
$N_m$. By (\ref{sp1}) we have 
 $$ 
 \lsp^{-1} \le
 \limsup (m/|X_{N_m}|)(1 + \limsup (K_m/m)).
 $$
 The first factor in the last expression is the inverse of the lim
inf speed in the graph $G_i$, for which we have the bound
(\ref{speedi}). Thus in order to show that $\lsp$ is greater than a
constant a.s. it suffices to find constants $b_\ell$ so that
 \begin{equation}
 \label{slowtime}
 \limsup (K_m/m) = \limsup (K_{m^2}/m^2) \le \sum_{\ell\ge 1} b_\ell
 <\infty \as
 \end{equation}
 The equality holds since $K_m$ is non-decreasing.

For each $\ell$, if $X_0=v$ is contained in a level $\ell$ country
$C$, then set $C_{\ell,0}:=C$, otherwise set
$C_{\ell,0}:=\emptyset$. Set $\tau_{\ell,0}:=0$, and for $k\ge 1$
define
 $$
 \tau_{\ell,k} = \min\{n\ge\tau_{\ell,k-1}+1: X_n \in C=:C_{\ell,k}
 \mbox{\ for some\  }C\in\ctry_\ell\setminus\{C_{\ell,k-1}\}\}.
 $$
 Also, for $k\ge 0$, define the time spent in the land between
stopping times:
 $$
 T_{\ell,k} = \#\{n: \tau_{\ell,k} \le n
 < \tau_{\ell,k+1},\ X_n\in C_{\ell,k}\cap A_i\}.
 $$
 We will use the rough bound
 $
 K_m \le \sum_{\ell=1}^{\infty} \sum_{k=0}^m T_{\ell,k}.
 $
 Since each vertex is contained in at most finitely many countries,
we have $T_{\ell,0}=0$ for all but finitely many
$\ell$. So for (\ref{slowtime}) it suffices to find summable
$b_\ell$ such that 
 \begin{equation}\label{ntb}
 \sum_{m\ge 1,\ell \ge 0}\pr(\sum_{k=1}^{m^2} T_{\ell,k}
 > b_\ell m^2)<\infty.
 \end{equation}

 Now fix $\ell$, and suppose that $X_0=v$ is on the inner vertex
boundary of a level $\ell$ country with land $R$. If $\rad(\ell)\ge
1$, then this means that $v$ is in the ocean and $\disti(v,
R)=\lfloor \rad(\ell) \rfloor$. Lemma \ref{time} with $i':=|R|^{-1}$
gives
 \begin{equation}
   \ev T_{\ell,0}^2 \le 8 (1-i^2)^{(\rad(\ell)-2)/2}|R|^{4.5} \le 8
   (1-i^2)^{a2^{\ell-1}/\ell^2-1}2^{4.5\ell} =:a^2_\ell.\label{al}
 \end{equation}
 If $\rad(\ell)<1$, then $v$ is contained in the land $R$, and this
bound still holds (although it is very rough) by Lemma \ref{time}
applied to a general starting point. By the Markov property, this
implies that for all $k\ge 1$, we have
 $
 \ev (T_{\ell,k}^2|{\cal F}(\tau_{\ell,k}))
 \le  a^2_\ell
 $, where ${\cal F}(\tau_{\ell,k})$ denotes the standard
 $\sigma$-field at the stopping time $\tau_{\ell,k}$, that
is the $\sigma$-field generated by information available up to
time $\tau_{\ell,k}$. Define
 $$
 S_{\ell,m}:=
 \sum_{k=1}^m \left(T_{\ell,k}
 - \ev \left(T_{\ell,k}|{\cal F}(\tau_{\ell,k})\right)\right)
 \ge  \sum_{k=1}^m (T_{\ell,k} - a_\ell) .
 $$
 Since $\{S_{\ell,m}\}_{m\ge 1}$ is a martingale, we can write
 $$
 \Var\ S_{\ell,m}
 = \sum_{k=1}^m
 \ev \Var \left(T_{\ell,k} | {\cal F}(\tau_{\ell,k})\right)
 \le m a^2_\ell.
 $$
 Therefore, if $b_\ell>a_\ell$, then Chebyshev's inequality gives
 $$
 \pr(\sum_{k=1}^{m} T_{\ell,k}> mb_\ell)
 \le \pr(S_{\ell,m}> m(b_\ell-a_\ell))
 \le \frac{m a_\ell^2}{m^2(b_\ell-a_\ell)^2}.
 $$
 Thus if we set, for example, $b_\ell:=(\ell+1) a_\ell$, then it is
clear from looking at the expression of $a_\ell$ that the conditions
of (\ref{ntb}) and (\ref{slowtime}) are satisfied. We thus have
proved that the speed is greater than a constant depending on $i$
and $\wmax$ only.

It remains to give a bound on the constant in terms of $i$
and $\wmax$. Since $b_\ell$ can be large when
$\ell$ is small, in order to get a reasonable bound, we need
to deal with  countries at or below some minimal level
$\ell_0$ separately.

The value of $\ell_0$ will be determined later, for now just assume
that $2^{-\ell_0}\le i$. Then by Lemma \ref{sink} the land of
countries of level up to $\ell_0$ is contained in $V\setminus
A_{2^{-\ell_0}}$. Let $K^*_m$ denote the time the random walk spends
in $V\setminus A_{2^{-\ell_0}}$ by time $N_m$, and let $K'_m$ denote
the time the walk spends in the land of countries of level greater
than $\ell_0$ by time $N_m$. We then have $N_m\le K^*_m+K'_m$.

Note that the sequence $\{K^*_m/|X_{N_m}|\}$ is a subsequence of
$\{m/|X_{N'_m}|\}$, where $N'_m$ is the time of the $m$th visit to
$G_{2^{-\ell_0}}$. By (\ref{speedi}), the lim sup of this sequence,
and thus the lim sup of the first one, is at most
$(2^{\ell_0})^2\log \wmax$. The bound (\ref{speedi}) on $\limsup
(m/|X_{N_m}|)$, (\ref{sp1}), and the bound (\ref{slowtime}) on
$\limsup (K'_\ell/m)$ from the first part of the proof imply
 \begin{eqnarray}\nonumber
 \lsp^{-1}
 &\le&
 \limsup (K^*_m/|X_{N_m}|) + \limsup (m/|X_{N_m}|)\limsup (K'_m/m)
 \\ &\le& \label{sp2}
 (2^{\ell_0})^2\log \wmax + i^{-2}\log\wmax
 \sum_{\ell>\ell_0}b_\ell.
 \end{eqnarray}
 We now want to choose an $\ell_0$ so that the last sum is small,
say each term $b_\ell$ is at most $2^{-\ell}$. From (\ref{al}),
since $\ell\ge 1$, we have
 \begin{equation}\label{2lbl}
 2^\ell b_\ell \le \exp(c\ell - c_2 \alpha^{-1} 2^\ell/\ell^2)
 \end{equation}
 with $\alpha:=\wmax i^{-3} \vee 2$. A simple computation shows
that there is a constant $c_3=c_3(c,c_2)\ge 1$ so that the right
hand side of (\ref{2lbl}) is at most $1$ if $\ell\ge\ell_0$, where
$\ell_0$ is chosen so that
 $$
 2^{ \ell_0} = c_3\alpha(\log_2 \alpha)^3.
 $$
 Using this choice of $\ell_0$, from (\ref{sp2}) we conclude that
 $$
 \lsp^{-1} \le c_3^2\alpha^2(\log_2 \alpha)^6 \log \wmax
 +i^{-2}\log\wmax
 \le c^{-1}\wmax^3\aec(G)^{-7}. \ \ \ \qed
 $$
 \def\qed{}
 \end{proof}

 \begin{eg}
 Consider the binary tree with edge weights $1$, and for each $n$
attach an extra vertex to each vertex at distance $n$ from the root
by an edge with weight $1/(n\log n)$. Add a self-loop to each new
vertex so that it will have weight 1. This graph clearly has
anchored expansion. The walk will visit infinitely many of these new
vertices by the Borel-Cantelli Lemma, and at each visit it has at
least constant probability to spend time at least $n \log n$ at the
vertex. Thus in this graph the speed is zero; this shows that in
Theorem \ref{main th} the bounded geometry condition cannot be left
out, nor replaced by bounds on the vertex weights.
 \end{eg}
 \begin{eg} 
 It follows from bounded geometry that positive transition
probabilities are bounded from below. This is another weaker
condition, but too weak for Theorem \ref{main th}. Define the {\bf
pipe} of length $n$ as the nearest neighbor graph on
$0,1,2,\ldots,n$. Consider the binary tree, and for every $n$ and
vertex at distance $n$, add a pipe of length $2k$ with edge weights
$1,2^{-1},\ldots,2^{-k+1}$,$2^{-k},2^{-k+1},\ldots,1$ with
$2^k\approx n\log n$. The argument of the previous example applies
again.
 \end{eg}

 \section{A heat kernel bound}
 \label{heat s}

This section contains the proof of Theorem \ref{heat th} and
examples showing that the bound there is sharp up to the constant
factor in the exponent.
 \begin{thm*}{\ref{heat th}} \heatth \end{thm*}

\begin{proof}
 Fix a vertex $x$, and let $a_n:=an^{1/3}$ for $a>0$ to be
determined later. Let $i:=\frac{2}{3}\aec(G)$, let $A_{i,n}\subset
A_i$ be the union of islands with volume at least $a_n$, and define
the {\bf territory} of such an island $C$ as the set of vertices
$v\in V$ with
 \begin{equation}
 \label{w2}
 \disti(v, C)\le a_n i/(4\wmax).
 \end{equation}
 Lemma \ref{connect} implies that there are only finitely many $n$
for which $x$ is contained in the territory of an island of
$A_{i,n}$. Consider large $n$ for which {\it (i)} this does not
happen and {\it (ii)} the right hand side of (\ref{w2}) is at least
$1$. Condition {\it (ii)} and the definition of $\disti$ ensures
that the inner vertex boundary of the territory of an island is a
subset of the ocean, $V\setminus A_i$.  Condition {\it (i)} implies
that $1/a_n<i$, and, as shown in Section \ref{geometry s},
$A_{1/a_n} \subset A_i$. By Lemma \ref{sink}, islands of $A_i$ with
volume less than $a_n$ do not intersect $A_{1/a_n}$, so
$A_{1/a_n}\subset A_{i,n}$, and we have $x\in V(G_{1/a_n})$.

Let $A'_{i,n}$ denote the union of islands of $A_{i,n}$ which are at
distance at most $n$ from $x$, and let $p'=p'(n)$ denote the transition
kernel of the walk on $G_{1/a_n}$. Note that $p^n(x,y)$ is the sum
of the probabilities of paths of length $n$ starting at $x$ and
ending at $y$. Each such path stays in $G_{1/a_n}$ or visits
$A'_{i,n}$. The total probability of the first kind of paths is at
most $p'^n(x,y)$ (regarded as $0$ if $y\in A_{1/a_n}$), so for all
$y$ we have
 $$
 p^n(x,y) \le p'^n(x,y)
 + \pr_x(\mbox{$\{X_k\}$ hits $A'_{i,n}$}).
 $$
The first term on the right satisfies
 \begin{eqnarray}\nonumber
 p'^n(x,y)
 &=&w(x)^{-1}(\one_x,P_{1/a_n}^n\one_y)
 \le \wmax^\hf\|P_{1/a_n}\|^n \\
 &\le& \wmax ^\hf (1-a_n^{-2})^{n \over 2}
 < \wmax^\hf \exp (-\sfrac{1}{2} na_n^{-2}).
 \label{Cheeger part}
 \end{eqnarray}
 The first inequality follows from Cauchy-Schwarz, the second from
Proposition \ref{Cheeger} and the fact that $G_{1/a_n}$ has Cheeger
constant at least $1/a_n$.

Suppose that there is a union $R_n$ of $n+1$ islands in $A'_{i,n}$
so that the territory of the first intersects the territory of all
the others. Then there is bridge structure $B$ interconnecting
$R_n\cup \{x\}$ with
 $$
 \#(B\cup\{x\}\setminus A_i)
 \le n + n(\sfrac{i}{2}\wmax^{-1} a_n-1) \le \wmax^{-1}\sfrac{i}{2}|R|.
 $$
 In the second expression, second term in parentheses is an upper
bound on the number of vertices in a bridge connecting two islands
with intersecting territories, and the first $n$ is an upper bound
for the number of vertices needed for the connection to $x$.  Lemma
\ref{connect} then implies that there are finitely many $n$ such
that $R_n$ exists. Thus for all large $n$, it is possible to
$n$-color islands in $A'_{i,n}$ so that islands of the same color
have disjoint territories. For such $n$, the probability of hitting
some island is bounded by $n$ times the maximal probability of
hitting some island of a given color.

 Suppose that the walk starts at a vertex $v$ on the inner vertex
boundary of the territory of an island $C \subset A'_{i,n}$. By
construction, this means that
 $$
 \disti (v, C) =
 \lfloor a_ni/(4\wmax)\rfloor \ge 1.
 $$
 Also note that another application of Lemma \ref{connect} shows
that for some $c$ and all large $n$, $A'_{i,n}$ cannot contain
islands with volume larger than $cn$. For such n, by Lemma
\ref{time} the probability of hitting $C$ is bounded by
 $$
 2w(x)^{-\hf}(1-i^2)^{(a_ni/(4\wmax)-2)/2} i^{-2} |\partial C|^\hf
 \le c(i,\wmax) n^{1/2}(1-i^2)^{(a_ni/(4\wmax)-2)/2}.
 $$
 In the first $n$ steps the walker has at most $n$ occasions to be
at the inner vertex boundary of some island of a given color. Thus
by the Markov property we get the bound
 \begin{equation}\label{Island part}
 \pr_x(\mbox{$\{X_k\}$ hits }A'_{i,n})
 \le n \cdot n \cdot c n^{1/2}\exp
 \left( \log(1-i^2)ia_n/( 8\wmax) \right).
 \end{equation}
 There exists an $a<(2\alpha)^{-1/2}$ so that the
exponents of (\ref{Cheeger part}) and (\ref{Island part}) are
at most $-cn^{-1/3}$
with $c>\alpha$. The statement of the theorem follows.
 \end{proof}

\begin{eg} 
 Let \walk\ be the nearest neighbor walk on the nonnegative integers
started at 0. There is a constant $a$ so that
$\pr(X_1,\ldots,X_{n^3} < n) > e^{-an}$ for all $n$. Consider the
binary tree with pipes of length $\ell_n$ attached to a vertex $v_n$
at distance $\ell_n$ from the root $o$ for some rapidly increasing
sequence $\{\ell_n\}$. This graph has anchored expansion. However,
consider the set of possible paths of length $\ell_n^3$ which start
and end at $o$. A subset of these start at $o$, travel on a shortest
path to the opposite end of the pipe starting at $v_n$, spend time
$\ell_n^3-4\ell_n$ in the pipe, and use the remaining time to return
to $o$. By the above, the probability measure of this set of paths
is at least $(1/4)^{4\ell_n}e^{-a\ell_n}$, and we get
$p^{\ell_n^3}(o,o) > e^{-c\ell_n}$. This shows that the conclusion
of Theorem \ref{heat th} is sharp up to the constant in the
exponent.
 \end{eg}

 \begin{eg}
 Chen and Peres (1998) showed that a supercritical Galton-Watson
tree has anchored expansion, so the above theorem gives the
$e^{-cn^{1/3}}$ heat kernel upper bound.  In this case, such bounds
are immediate from results of Piau. For Galton-Watson trees where
the probability of non-branching (zero or one offspring) is positive
this bound is easily seen to be sharp up to the constant in the
exponent (see Piau 1998).
 \end{eg}

\noindent {\bf Acknowledgments.} 
 The author thanks Itai Benjamini, Russell Lyons, Yuval Peres and
Oded Schramm for illuminating discussions and comments on previous
versions. Special thanks are due to Russell Lyons and the referee
for their help in eliminating many mistakes and making this paper
more readable.

\section*{References}
\newcommand{\andd}{and\ }
\def\labelenumi{[\arabic{enumi}]}
\begin{enumerate}
\item
Benjamini, I., Lyons, R. \andd Schramm, O. (1999)
Percolation perturbations in potential theory and random walks,
In Picardello, M. and Woess, W., editors,
{\it  Random Walks and Discrete Potential Theory (Cortona, 1997)}, 56--84.
Sympos. Math. Cambridge Univ. Press, Cambridge.
\item
Cheeger, J.
(1970) A lower bound for the lowest eigenvalue of the Laplacian,
{\it Problems in Analysis (Sympos. Salomon Bochner, Princeton Univ., 1969),}
 195--199. Princeton Univ. Press, Princeton, NJ.
\item
Chen, D. \andd Peres, Y. (1998) Anchored expansion, percolation and speed,
{\it Preprint.}
\item
Dodziuk, J. (1984) Difference equations, isoperimetric inequalities,
and transience of certain random walks, {\it Trans. Amer. Math. Soc.}
{\bf 284}, 787-794.
\item
H\"aggstr\"om, O., Schonmann, R. and Steif, J. (1999)
The Ising model on diluted graphs and strong amenability,
{\it Ann. Probab.} To appear.
\item
Hebisch W. and Saloff-Coste L. (1993) Gaussian estimates for Markov
chains and random walks on groups, {\it Ann. Probab.} {\bf 21}, 673--709.
\item
Lyons, R. \andd Peres, Y. (1998) {\it Probability on Trees and Networks}
(a book),
Cambridge University Press, in progress. Current version published on
the web at \\ {\tt
http://php.indiana.edu/\~{}rdlyons}.
\item
Mohar, B. (1988) Isoperimetric inequalities, growth, and the spectrum
of graphs, Linear Algebra Appl. {\bf 103}, 119--131.
\item
Thomassen, C. (1992) Isoperimetric inequalities and transient
random walks on graphs, {\it Ann. Probab.} {\bf 20}, 1592-1600.
\item
Piau, D. (1998) Functional central limit theorem for a random walk in a
random environment. {\it Ann. Probab.} {\bf 26}, 1016--1040.
\item
Varopoulos, N. Th. (1991) Groups of superpolynomial growth, in {\it Harmonic
analysis (Sendai, 1990)}, 194--200, Springer, Tokyo.

\end{enumerate}
\noindent Department Of Mathematics, Massachusetts Institute of Technology, 
Cambridge, MA 02139, USA
 \\ {\tt balint@math.mit.edu}

\end{document}